\documentclass{amsart}
\usepackage{graphicx}
\usepackage{amssymb}
\usepackage{amsfonts}
\swapnumbers
\newtheorem{theorem}{Theorem}[section]

\theoremstyle{definition}

\newtheorem{remark}[theorem]{Remark}

\numberwithin{equation}{section}
 \theoremstyle{plain}
 
 \numberwithin{equation}{section} %% Comment out for sequentially-numbered
 \numberwithin{figure}{section} %% Comment out for sequentially-numbered
 \theoremstyle{plain}
 \theoremstyle{plain}
 \theoremstyle{remark}
 \newtheorem*{acknowledgement*}{Acknowledgement}
\newcommand{\Om}{{\Omega}}

\newcommand{\ve}{{\varepsilon}}
\newcommand{\del}{{\delta}}
\newcommand{\Del}{{\Delta}}

\newcommand{\Gam}{{\Gamma}}

\newcommand{\Sig}{{\Sigma}}
\newcommand{\sig}{{\sigma}}
\newcommand{\al}{{\alpha}}
\newcommand{\be}{{\beta}}

\newcommand{\la}{{\lambda}}

\newcommand{\bbC}{{\mathbb C}}

\newcommand{\bbR}{{\mathbb R}}

\newcommand{\bbZ}{{\mathbb Z}}

\begin{document}
\title[]{An Erd\" os-R\' enyi law for nonconventional sums}%
 \vskip 0.1cm
 \author{Yuri Kifer\\
\vskip 0.1cm
 Institute  of Mathematics\\
Hebrew University\\
Jerusalem, Israel}%
\address{
Institute of Mathematics, The Hebrew University, Jerusalem 91904, Israel}
\email{ kifer@math.huji.ac.il}%

\thanks{ }
\dedicatory{  }
 \date{\today}
\begin{abstract}\noindent
We obtain the Erd\" os-R\' enyi type law of large numbers for
"nonconventional" sums of the form $S_n=\sum_{m=1}^n
F(X_m,X_{2m},...,X_{\ell m})$ where $X_1,X_2,...$ is a sequence of i.i.d.
random variables and $F$ is a bounded Borel function. The proof relies
on nonconventional large deviations obtained in \cite{KV}.
\end{abstract}
%\footnotetext[1]{}
\maketitle
\markboth{Y. Kifer }{Erd\" os-R\' enyi law}
\renewcommand{\theequation}{\arabic{section}.\arabic{equation}}
\pagenumbering{arabic}

\renewcommand{\theequation}{\arabic{section}.\arabic{equation}}
\pagenumbering{arabic}

\section{Introduction}\label{sec1}\setcounter{equation}{0}

Let $X_1,X_2,...$ be a sequence of independent identically distributed
(i.i.d.) random variables such that $EX_1=0$ and the moment generating
function $\phi(t)=Ee^{tX_1}$ exists. Denote by $I$ the Legendre transform
of $\ln\phi$ and set $S_n=\sum_{m=1}^nX_m$ for $n\geq 1$ and $S_0=0$.
The Erd\" os-R\' enyi law of large
 numbers from \cite{ER} says that with probability one
\begin{equation}\label{1.1}
I(\al)\lim_{n\to\infty}\max_{0\leq m\leq n-[\frac {\ln n}{I(\al)}]}\frac
{S_{m+[\frac {\ln n}{I(\al)}]}-S_m}{\ln n}=\al
\end{equation}
for all $\al>0$ in some neighborhood of zero.

The nonconventional limit theorems initiated in \cite{Ki} and partially
motivated by nonconventional ergodic theorems study asymptotic behaviors
of sums of the form $S_n=\sum_{m=1}^nF(X_m,X_{2m},...,X_{\ell m})$ and
more general ones where $F$ is a Borel function. In this paper we will
obtain an Erd\" os-R\' enyi law similar to (\ref{1.1}) for such sums where
$X_1,X_2,...$ is again a sequence of i.i.d. random variables and $F$ is a
bounded Borel function. Observe that summands in nonconventional sums are
long range dependent so this result cannot be derived directly from existing
literature. On the other hand, as most proofs of the Erd\" os-R\' enyi
law we will rely on large deviations which in the nonconventional setup
were obtained in \cite{KV}.

\section{Preliminaries and main results}\label{sec2}\setcounter{equation}{0}

Let $X_1,X_2,...$ be a sequence of i.i.d. random variables and $F$ be a
bounded Borel function on $\bbR^\ell$ such that
\begin{equation}\label{2.1}
\bar F=EF(X_1,X_2,...,X_\ell)=0\quad\mbox{and}\quad\sig^2=EF^2(X_1,X_2,...,
X_\ell)>0.
\end{equation}
The first condition in (\ref{2.1}) is not a restriction since we always can
consider $F-\bar F$ in place of $F$ and the second condition there means that
$F$ is not a constant almost surely (a.s.) with respect to the $\ell$-product
measure
$\mu^{(\ell)}=\mu\times\mu\times\cdots\times\mu$ on $\bbR^\ell$ where $\mu$
is the distribution of $X_1$. Set also $M=\| F\|_\infty$ and $M_+=
\| F_+\|_\infty$ where $F_+(x_1,...,x_\ell)=\max(0,F(x_1,...,x_\ell))$ and
the $L^\infty$ norm on $\bbR^\ell$ is considered with respect to the measure
 $\mu^{(\ell)}$.
Introduce the moment generating function $\phi(t)=E\exp(tF(X_1,X_2,...,X_\ell))$
and its Legendre transform
\begin{equation}\label{2.2}
I(\al)=sup_t(t\al-\ln\phi(t)).
\end{equation}
\begin{theorem}\label{thm2.1}
 With $I$ given by (\ref{2.2}) the Erd\" os-R\' enyi law
(\ref{1.1}) holds true also for the nonconventional sums $S_n=\sum_{m=1}^n
F(X_m,X_{2m},...,X_{\ell m})$ for all $\al\in(0,M_+)$ where we set also $S_0=0$.
\end{theorem}

Our proof of Theorem \ref{thm2.1} will follow the scheme of \cite{DK} but
we will rely also on nonconventional large deviations results from \cite{KV}.
As in some books and many papers on large deviations we did not address
explicitly
in \cite{KV} the crucial question when the rate function of large deviations
is positive without which the large deviations principle is meaningless since
it does not lead to any nontrivial estimates for the domains were the rate
function is zero. We will rely on the following theorem which specifies
further the results of \cite{KV} and actually provides more information than
we need for the proof of Theorem \ref{thm2.1}.

\begin{theorem}\label{thm2.2} The limit
\begin{equation}\label{2.3}
Q(\la F)=\lim_{N\to\infty}\frac 1N\ln E\exp(S_N(\la F))
\end{equation}
exists where $Q(\la F)$ is a $\bbC^\infty$ function of $\la$ with bounded
derivatives and $S_n,\, n\geq 1$ are nonconventional sums from Theorem
\ref{thm2.1}. The Legendre transform of $Q$,
\begin{equation}\label{2.4}
J(u)=\sup_\la(\la u-Q(\la F))
\end{equation}
is a nonnegative, convex, lower semi-continuous function such that
$J(u)=0$ if and only if $u=0$ and $J(u)$ is strictly increasing for $u\geq 0$
( writing for convenience $\infty>\infty$)
 while it is strictly decreasing for $u\leq 0$. In addition,
 if $M_+=\| F_+\|_\infty>0$ ($M_-=\| F-F_+\|_\infty>0$) then there exists
 $L_+>0$ ($L_->0$) such that $J(u)<\infty$ when $u\in[0,L_+)$ ($u\in
 (-L_-,0]$) and $J(u)=\infty$ when $u>L_+$ ($u<-L_-$).
  Furthermore, the sums
   $S_n,\, n\geq 1$ satify the large deviations principle in the form
\begin{equation}\label{2.5}
\limsup_{N\to\infty}\frac 1N\ln P\{\frac 1NS_N\in K\}\leq -\inf_{u\in K}J(u)
\end{equation}
for any closed set $K\subset\bbR$ while for any open set $U\subset\bbR$,
\begin{equation}\label{2.6}
\liminf_{N\to\infty}\frac 1N\ln P\{\frac 1NS_N\in U\}\geq -\inf_{u\in U}J(u).
\end{equation}
\end{theorem}

\begin{remark}\label{rem2.3}
Theorem \ref{thm2.1} shows that the Erd\" os-R\' enyi law for
nonconventional sums has the same form as for sums of i.i.d. random variables
having the same distribution as $F(X_1,X_2,...,X_\ell)$. This is similar
to the nonconventional strong law of large numbers proved in \cite{Ki1}.
On the other hand, the nonconventional central limit theorem and the
nonconventional large deviations estimates are somewhat different from the
corresponding results for sums of i.i.d. random variables. In particular,
it is shown in \cite{KV0} that the nonconventional functional central limit
theorem may yield in the limit a process with dependent increments while
concerning large deviations it follows from \cite{KV} that the rate functions
 $I$ and $J$ above are, in general, different.
\end{remark}

\section{Proof of Theorem \ref{thm2.1}}\label{sec3}\setcounter{equation}{0}

Let $Y_1,Y_2,...$ be a sequence of i.i.d. random variables which have
the same distribution as $F(X_1,X_2,...,X_\ell)$ and set $\Sig_n=
\sum_{m=1}^nY_m$. We will need the classical Cram\' er large deviation
estimates in the form (see, for instance, Section 2.2 in \cite{DZ}),
\begin{equation}\label{3.1}
\limsup_{N\to\infty}\frac 1N\ln P\{\frac 1N\Sig_N\in K\}\leq -\inf_{u\in K}I(u)
\end{equation}
for any closed set $K\subset\bbR$ while for any open set $U\subset\bbR$,
\begin{equation}\label{3.2}
\liminf_{N\to\infty}\frac 1N\ln P\{\frac 1N\Sig_N\in U\}\geq -\inf_{u\in U}I(u)
\end{equation}
where $I$ is given by (\ref{2.2}).

It is essential to observe that $I(\al)>0$ ($I(\al)=\infty$ is possible)
unless $\al=0$ which is well known and follows, in particular, from Theorem
II.6.3 in \cite{El} (which relies on general convex analysis results)
 but it has also a simple direct explanation in our case. Indeed, since
 $\ln\phi(0)=(\ln\phi(t))'_{t=0}=0$ then $\ln\phi(t)=o(t)$ for small $t$. Hence,
 if $\al\ne 0$ then $t\al>\ln\phi(t)$ either for small positive or for small
 negative $t$, and so in view of (\ref{2.2}), $I(\al)=0$ only when $\al=0$ and
 otherwise $I(\al)$ is positive. By (\ref{2.1}) and the Jensen inequality
 $\ln\phi(t)\geq tEF(X_1,...,X_\ell)=0$, and so (see Lemma 2.2.5 in
 \cite{DZ}),
 \begin{equation}\label{3.2+}
 I(\al)=sup_{t\geq 0}(t\al-\ln\phi(t))\,\,\mbox{if}\,\,\al\geq 0\,\,\mbox{and}
 \,\, I(\al)=sup_{t\leq 0}(t\al-\ln\phi(t))\,\,\mbox{if}\,\,\al\leq 0.
 \end{equation}
 For each $\al>0$ there exists a sequence $t_n\to t_0$ as $n\to\infty$
 such that $I(\al)=\lim_{n\to\infty}(t_n\al-\ln\phi(t_n))$ where $t_0>0$
 ($t_0=\infty$ is possible) since by above $I(\al)>0$. Therefore, for
 any $\Del>0$,
 \[
 I(\al+\Del)\geq\lim_{n\to\infty}(t_n(\al+\Del)-\ln\phi(t_n))=I(\al)+t_0\Del
 \]
 which means that $I(\al)$ is strictly increasing for $\al\geq 0$. Similarly,
 $I(\al)$ is strictly decreasing for $\al\leq 0$. Observe that, in fact, for
 any $\ve>0$,
 \[
 e^{tM_+}+1\geq\phi(t)\geq P\{ F(X_1,...,X_\ell)\geq M_+-\ve\}e^{t(M_+-\ve)}
 \,\,\mbox{if}\,\, t\geq 0\,\,\mbox{and}
 \]
 \[
 e^{-tM_-}+1\geq\phi(t)\geq P\{ -F(X_1,...,X_\ell)\geq M_--\ve\}e^{-t(M_--\ve)}
 \,\,\mbox{if}\,\, t\leq 0
 \]
 where $M_-=\| F-F_+\|_\infty$. This together with (\ref{2.2}) yields
 that $I(\al)<\infty$
  if $-M_-<\al<M_+$ while $I(\al)=\infty$ if $\al>M_+$ or $\al<-M_-$.
 Similar arguments relying on explicit formulas from \cite{KV} yield Theorem
  \ref{thm2.2} but for now we will take it for granted in order to prove
   Theorem \ref{thm2.1}.

 Fix $\al\in(0,M)$ and let $b_n=[\ln n/I(\al)]$. Choose $\ve>0$ and define the
 event
 \[
 A_n(\ve)=\{\max_{0\leq m\leq n-b_n}(S_{m+b_n}-S_m)\geq(\al+\ve)b_n\}.
 \]
 Then
\begin{eqnarray}\label{3.3}
&P(A_n(\ve))=P\big\{\bigcup_{0\leq m\leq n-b_n}\{S_{m+b_n}-S_m\geq(\al+\ve)
b_n\}\big\}\\
&\leq\displaystyle\sum_{m=0}^{n-b_n}P\{ S_{m+b_n}-S_m\geq(\al+\ve)b_n\}.
\nonumber\end{eqnarray}

Observe that when $m> (\ell-1)b_n$ then
\[
 S_{m+b_n}-S_m=\sum_{k=m+1}^{m+b_n}F(X_k,X_{2k},...,X_{\ell k})
 \]
 is the sum of i.i.d. random variables having the same distribution as
 $F(X_1,X_2,...,X_\ell)$. Indeed, if $\ell=1$ this is clear and if $\ell>1$
 then the equality $ik=j\tilde k$ is impossible for integers $(\ell-1)b_n\leq
 m<\tilde k<k\leq m+b_n$ and $1\leq i<j\leq\ell$ since $j/i\geq 1+
 (\ell-1)^{-1}$ while $k/\tilde k<1+b_n/m\leq 1+(\ell -1)^{-1}$.
 Hence, we can use Cram\' er's upper large deviations bound
 (\ref{3.1}) to conclude that for any $m\geq (\ell-1)b_n$,
 \begin{equation}\label{3.4}
 P\{ S_{m+b_n}-S_m\geq(\al+\ve)b_n\}\leq\exp(-b_n(I(\al+\ve)-\del))\leq
 \exp(-b_n(I(\al)+\del))
 \end{equation}
 where $0<\del<\frac 12(I(\al+\ve)-I(\al))$, $n\geq n(\del)$ is large enough
  and we use the fact that $I(\be)$ is strictly increasing when $\be\geq 0$.

 On the other hand, if $m\leq (\ell-1)b_n$ and $\ell>1$ then we write
 \begin{eqnarray}\label{3.5}
 &P\{ S_{m+b_n}-S_m\geq(\al+\ve)b_n\}\leq P\{ S_{m+b_n}\geq\frac
 12(\al+\ve)b_n\}\\
 &+P\{ -S_m\geq\frac 12(\al+\ve)b_n\}\leq P\{\frac 1{m+b_n}S_{m+b_n}\geq
 \frac 1{2\ell}(\al+\ve)\}\nonumber\\
 &+ P\{ -\frac 1mS_m\geq\frac 1{2m}(\al+\ve)b_n\}.\nonumber
 \end{eqnarray}
 Applying the upper nonconventional large deviations bound (\ref{2.5})
 we obtain for $m\leq (\ell-1)b_n$ that
 \begin{eqnarray}\label{3.6}
 &P\{ \frac 1{m+b_n}S_{m+b_n}\geq\frac 1{2\ell}(\al+\ve)\}\\
 &\leq\exp(-(m+b_n)(J(\frac 1{2\ell}(\al+\ve))-\del))\leq
 \exp(-\frac 12b_nJ(\frac \al{2\ell}))\nonumber
 \end{eqnarray}
 where $0<\del<J(\frac 1{2\ell}(\al+\ve))-J(\frac \al{2\ell})$, $n\geq n(\del)$
 is large enough and we use that $J(\be)$ is strictly increasing when
 $\be\geq 0$.

 Since $|S_m|\leq mM$ a.s. then
 \begin{equation}\label{3.7}
 P\{ -\frac 1mS_m\geq\frac 1{2m}(\al+\ve)b_n\}=0\,\,\,\mbox{if}\,\,\,
 m<\frac \al{2M}b_n.
 \end{equation}
 Now assume that $(\ell-1)b_n\geq m\geq\frac \al{2M}b_n$ and $\ell>1$.
 Observe that
 $-S_m=\sum_{k=1}^m(-F(X_k,X_{2k},...,X_{\ell k}))$, and so we can consider
 nonconventional large deviations estimates of Theorem \ref{thm2.2} for the
 case where $F$ is replaced by $-F$ with a corresponding rate function
 $\hat J$ having the same properties as $J$. Then we obtain
 \begin{eqnarray}\label{3.8}
 &P\{ -\frac 1mS_m\geq\frac 1{2m}(\al+\ve)b_n\}\\
 &\leq P\{ -\frac 1mS_m\geq\frac 1{2(\ell-1)}(\al+\ve)\}\nonumber\\
 &\leq\exp(-m(\hat J(\frac 1{2(\ell-1)}(\al+\ve))-\del))\leq
 \exp(-b_n\frac \al{2M}\hat J(\frac \al{2\ell}))\nonumber
 \end{eqnarray}
 where $0<\del<\hat J(\frac 1{2(\ell-1)}(\al+\ve))-\hat J(\frac \al{2\ell})$,
  $n\geq n(\del)$ is large enough and we use that $\hat J(\be)$ is strictly
  increasing when $\be\geq 0$ (of course, if $\hat J(\frac 1{2(\ell-1)}
  (\al+\ve))=\infty$ then any $\del$ will do).

 Set $c=c_\al=\frac 12\min(J(\frac \al{2\ell}),\frac \al{M}
 \hat J(\frac \al{2\ell}))$
 which is a positive number. Then it follows from (\ref{3.3})--(\ref{3.8})
 that for $\del$ satisfying (\ref{3.4}) and for $n$ large enough,
 \begin{eqnarray}\label{3.9}
 &P(A_n(\ve))\leq n\exp(-b_n(I(\al)+\del))+\ell b_n\exp(-c b_n)\\
 &\leq n\exp\big(-(\frac {\ln n}{I(\al)}-1)(I(\al)+\del)\big)
 +\ell(\frac {\ln n}{I(\al)}+1)\exp(-c(\frac {\ln n}{I(\al)}-1))
 \nonumber\\
 &=e^{I(\al)+\del}n^{-\frac \del{I(\al)}}+\ell(\frac {\ln n}{I(\al)}+1)e^{c}
 n^{-\frac {c}{I(\al)}}.
  \nonumber\end{eqnarray}
  Now let $d>I(\al)\max(\del^{-1},c^{-1})$. Then
  \[
  \sum_{n=1}^\infty(n^{-\frac {d\del}{I(\al)}}+n^{-\frac {dc}{I(\al)}}\ln n)
  <\infty
  \]
 which together with (\ref{3.9}) and the Borel-Cantelli lemma yields that
 with probability one $A_{n^d}(\ve)$ occurs only finitely often. Hence,
 setting $a_n=b_{n^d}$ we obtain
 \begin{equation}\label{3.10}
 \limsup_{n\to\infty}\max_{0\leq m\leq n^d-a_{n}}\frac {S_{m+a_{n}}-S_m}
 {a_{n}}\leq\al+\ve.
 \end{equation}
 Since for $n^d<r\leq (n+1)^d$ large enough the difference $a_{n}-b_r$
 is bounded by 1 then it follows that
 \begin{eqnarray}\label{3.11}
 &\limsup_{r\to\infty}\max_{0\leq m\leq r-b_{r}}\frac {S_{m+b_{r}}-S_m}
 {b_{r}}\\
 &\leq\limsup_{n\to\infty}\max_{0\leq m\leq (n+1)^d-a_{n+1}}\frac
 {S_{m+a_{n+1}}-S_m+M}{a_{n}}\leq\al+\ve.\nonumber
 \end{eqnarray}

 In order to derive the lower bound choose $\ve>0$ so that $\al-\ve>0$
 and define
 \[
 B_n(\ve)=\{\max_{0\leq m\leq n-b_n}(S_{m+b_n}-S_m)\leq b_n(\al-\ve)\}.
 \]
 Let $C_m=\{ S_{m+b_n}-S_m\leq b_n(\al-\ve)\}$. Then
 \begin{equation}\label{3.12}
 P(B_n(\ve))=P(\bigcap_{0\leq m\leq n-b_n}C_m)\leq P(\bigcap_{(1-\ell^{-1})
 n\leq m\leq n-b_n}C_m).
 \end{equation}
 Observe that when $n-b_n\geq m,\tilde m\geq (1-\ell^{-1})n$ then
 $\frac m{\tilde m}<\frac \ell{\ell-1}$ if $m>\tilde m$, and so the
 equality $im=j\tilde m$ for integers $n\geq m>\tilde m\geq (1-\ell^{-1})n$
 and $\ell\geq j>i\geq 1$ is impossible since then $\min\frac ji=
 \frac \ell{\ell-1}$. Hence, all $F(X_k,X_{2k},...,X_{\ell k})$,
 $(1-\ell^{-1})n\leq k\leq n-b_n$ are independent, and so all events
 $C_{mb_n},\, (1-\ell^{-1})n\leq mb_n\leq n-b_n$ are independent for different m. 
 Hence by
 (\ref{3.12}),
 \begin{equation}\label{3.13}
 P(B_n(\ve))\leq\prod_{m:\,(1-\ell^{-1})n\leq mb_n\leq n-b_n}P(C_{mb_n}).
 \end{equation}

 Taking into account that $S_{m+b_n}-S_m$ is a sum of i.i.d. random variables
 having the same distribution as $F(X_1,X_2,...,X_\ell)$ when
 $(1-\ell^{-1})n\leq m\leq n-b_n$ we obtain by Cram\' er's lower large
 deviations bound (\ref{3.2}) that
 \begin{equation}\label{3.14}
 P(\Om\setminus C_m)\geq\exp(-b_n(I(\al-\ve)+\del))\geq\exp(-b_nI(\al)(1-\del))
 \geq\exp(-(1-\del)\ln n)
 \end{equation}
 where we choose $\del>0$ so small that $(I(\al-\ve)+\del)/I(\al)<1-\del$ which
 is possible since $I(\be)$ is strictly increasing for $\be\geq 0$. Hence,
 if $n$ is sufficiently large,
 \begin{eqnarray}\label{3.15}
 &P(B_n(\ve))\leq (1-\exp(-(1-\del)\ln n)^{\frac n{\ell b_n}-2}=
 (1-n^{-(1-\del)})^{\frac n{\ell b_n}-2}\\
 &=\big((1-n^{-(1-\del)})^{n^{1-\del}}\big)^{\frac {n^\del}{\ell b_n}-1}=
 O(\exp(-n^{\del/2}).\nonumber
 \end{eqnarray}
 It follows that
 \[
 \sum_{n=1}^\infty P(B_n(\ve))<\infty
 \]
 and by the Borel-Cantelli lemma with probability one $B_n(\ve)$ occurs
 only finitely often which implies that
 \begin{equation}\label{3.16}
 \liminf_{n\to\infty}\max_{0\leq m\leq n-b_n}\frac {S_{m+b_n}-S_m}{b_n}\geq
 \al+\ve.
 \end{equation}
 Since $\ve>0$ can be chosen arbitrarily small we obtain the assertion of
 Theorem \ref{thm2.1} from (\ref{3.11}) and (\ref{3.16}).   \qed

 \section{Proof of Theorem \ref{thm2.2}}\label{sec4}\setcounter{equation}{0}

Theorem \ref{thm2.2} mostly follows from the results of \cite{KV} together
with Theorem II.6.3 from \cite{El} but for reader's convenience we will give a
direct argument here. First, we recall relevant notations and formulas from
\cite{KV}. Let $r_1,...,r_m\geq 2$ be all primes not exceeding $\ell$. Set
$A_n=\{ a\leq n:\, a$ is coprime with $r_1,...,r_m\}$ and $B_n(a)=\{ b\leq n:\,
b=ar_1^{d_1}r_2^{d_2}\cdots r_m^{d_m}$ for some nonnegative integers
$d_1,...,d_m\}$. For any function $V$ on $\bbR^\ell$ we write
\[
S_N(V)=\sum_{a\in A_N}S_{N,a}(V)\,\,\mbox{where}\,\, S_{N,a}(V)=\sum_{b\in
B_N(a)}V(X_b,X_{2b},...,X_{\ell b})
\]
observing that $S_n$ from Theorem \ref{thm2.2} equals $S_n(F)$ here.

The existence of the limit (\ref{2.3}) was proved in \cite{KV}.
 Recall, that convexity and lower semi-continuity of the Legendre transform
 $J(u)$ of $Q(\la F)$ follows
from (\ref{2.3}) and (\ref{2.4}) automatically (see Theorem II.6.1 in
\cite{El}). Observe that by (\ref{2.1}) and the Jensen
 inequality $Q(\la F)\geq 0$ and since $Q(\la F)\leq |\la| M$ then
 $J(u)=\infty$ when $u>M$.
Note that Theorem 2.7 in \cite{KV} is formulated for continuous
functions but, in fact, only boundedness of functions is used in the proof
so we can apply it to our setup where $\| F\|_\infty=M<\infty$.

In order to exhibit an explicit formula for $Q(\la F)$ obtained in \cite{KV}
introduce
\[
D(\rho)=\{ n=(n_1,...,n_m)\in\bbZ^m:\, n_1,...,n_m\geq 0,\,\,\mbox{and}\,\,
\sum_{i=1}^mn_i\ln r_i\leq\rho\}
\]
and observe that $D(\ln(N/a))|=|B_N(a)|=|B_{N/a}(1)|$ where $|\Gam|$ denotes
the cardinality of a finite set $\Gam$. Set
\[
\rho_{\mbox{\tiny min}}(l)=\inf\{\rho\geq 0:\, |D(\rho)|=l\}\,\,\mbox{and}\,\,
\rho_{\mbox{\tiny max}}(l)=\sup\{\rho\geq 0:\, |D(\rho)|=l\}.
\]
It was shown in \cite{KV} that for any $l\geq 1$,
\begin{equation}\label{4.3}
\rho_{\mbox{\tiny max}}(l)>\rho_{\mbox{\tiny min}}(l)\geq (l^{1/m}-1)\ln 2.
\end{equation}

Set
\[
Z_{n,a}(\la F)=E\exp S_{N,a}(\la F).
\]
As it was explained in \cite{KV} the distribution of $S_{N,a}(\la F)$ depends
only on $|B_n(a)|$ (in addition to $\la F$, of course), and so $Z_{n,a}(\la F)$
 is determined by $|B_n(a)|$. Hence,
we can set $R_l(\la F)=Z_{n,a}(\la F)$ provided $|B_n(a)|=l$. Now we can write
the formula for $Q$ obtained in \cite{KV},
\begin{equation}\label{4.4}
Q(\la F)=r\sum_{l=1}^\infty(e^{-\rho_{\mbox{\tiny min}}(l)}-e^{-
\rho_{\mbox{\tiny max}}(l)})\ln R_l(\la F)
\end{equation}
where
\[
r=\prod_{k=1}^m(1-\frac 1{r_k})=1+\sum_{k=1}^m(-1)^k\sum_{i_1<i_2<...<i_k\leq m}
\prod_{j=1}^k\frac 1{r_{i_j}}.
\]
The series in (\ref{4.4}) converges absolutely in view of (\ref{4.3}) taking
into account that $\ln R_l(\la F)\leq lM|\la|$. By (\ref{2.1}) and the Jensen
inequality we have also that $\ln R_l(\la F)\geq 0$.

Now observe that for any $k\geq 1$,
\begin{equation}\label{4.5}
\big\vert\frac {d^kR_l(\la F)}{d\la^k}\big\vert\leq l^kM^kR_l(\la F),
\end{equation}
and so
\begin{equation}\label{4.6}
\big\vert\frac {d^k\ln R_l(\la F)}{d\la^k}\big\vert\leq C_kl^kM^k
\end{equation}
for some $C_k>0$ depending only on $k$. It follows that $Q(\la F)$ is
$\bbC^\infty$ in $\la$ and
\begin{equation}\label{4.7}
\big\vert\frac {d^kQ(\la F)}{d\la^k}\big\vert\leq\hat C_k
\end{equation}
where
\[
\hat C_k=C_krM^k\sum_{l=1}^\infty(e^{-\rho_{\mbox{\tiny min}}(l)}-
e^{-\rho_{\mbox{\tiny max}}(l)})l^k
\]
and the latter series converges absolutely in view of (\ref{4.3}).
Note that existence of the first derivative of $Q(\la F)$ in $\la$
already yields the large deviations bounds (\ref{2.3}) and (\ref{2.4})
(see Theorem II.6.1 in \cite{El}).

Now observe that in view of (\ref{2.1}),
\[
\frac {d\ln R_l(\la F)}{d\la}\big\vert_{\la=0}=\frac {dR_l(\la F)}{d\la}
\big\vert_{\la=0}=0,
\]
and so
\begin{equation}\label{4.8}
\frac {dQ(\la F)}{d\la}\big\vert_{\la=0}=0.
\end{equation}
From Theorem II.6.3 in \cite{El} it follows that (\ref{2.1}), (\ref{2.3}),
(\ref{2.4}) and
(\ref{4.8}) yield already that $J(u)$ attains its infimum at the unique point 0
and it is positive when $|u|>0$.
As in Section \ref{sec3} the direct argument proceeds as follows. Since $Q(0)=0$
then  (\ref{4.8}) implies that $Q(\la F)=o(\la)$ for small $\la$, and so
$|\la u|>Q(\la F)$ when $|\la|$ is small which together with (\ref{2.4})
yields the assertion above.

 Taking into account that $Q(\la F)\geq 0$ we see that
 \begin{equation}\label{4.9}
 J(u)=sup_{\la\geq 0}(\la u-Q(\la F))\,\,\mbox{if}\,\, u\geq 0\,\,\mbox{and}
 \,\, J(u)=sup_{\la\leq 0}(\la u-Q(\la F))\,\,\mbox{if}\,\, u\leq 0.
 \end{equation}
 Similarly to Section \ref{sec3} we argue that for each $u>0$ there exists
  a sequence $\la_n\to \la_0$ as $n\to\infty$
 such that $J(u)=\lim_{n\to\infty}(\la_n u-Q(\la F))$ where $\la_0>0$
 ($t_0=\infty$ is possible) since by above $J(u)>0$. Therefore, for
 any $\Del>0$,
 \[
 J(u+\Del)\geq\lim_{n\to\infty}(\la_n(u+\Del)-Q(\la_nF))=J(u)+\la_0\Del
 \]
 which means that $J(u)$ is strictly increasing for $u\geq 0$. Similarly,
 $J(u)$ is strictly decreasing for $\leq 0$.

 Observe that by Jensen's inequality $\ln R_l(\la F)\geq 0$ for all $l\geq 1$,
and so all terms of the series in (\ref{4.4}) are nonnegative. Hence,
for any $\la>0$ and $\ve\in(0,M_+)$,
\[
Q(\la F)\geq K\ln R_1(\la F)\geq L(\ve)\la>0
\]
where $K=r(e^{-\rho_{\mbox{\tiny min}}(1)}-e^{-\rho_{\mbox{\tiny max}}(1)})$
and $L(\ve)=KP\{ F(X_1,...,X_\ell)\geq M_+-\ve\}(M_+-\ve)$. Then, clearly,
 $J(u)<\infty$ for all $u\in[0,L(\ve))$. If $M_+>0$ then
 by the monotonicity property of $J$
 obtained above we conclude that there exists $L_+>0$ such that $J(u)<\infty$
 for $u\in[0,L_+)$ while $J(u)=\infty$ for $u>L_+$. Similarly, if $M_->0$
 then there exists
 $L_->0$ such that $J(u)<\infty$
 for $u\in(-L_-,0]$ while $J(u)=\infty$ for $u<-L_-$, completing the proof
 of Theorem \ref{2.2}.    \qed

\end{document}